# ON CONVERGENCE ALMOST EVERYWHERE OF MULTIPLE FOURIER INTEGRALS
(Mengenai Penumpuan Hampir di Mana-Mana bagi Kamiran Fourier Berganda)

ANVARJON AHMEDOV, NORASHIKIN ABDUL AZIZ & MOHD NORIZNAN MOHTAR

## ABSTRACT

In this paper we investigate the principle of the generalised localisation for spectral expansions of the polyharmonic operator, which coincides with the multiple Fourier integrals summed over the domains corresponding to the surface levels of the polyharmonic polynomials. It is proved that the partial sums of the multiple Fourier integrals of a function $f \in L_2(R^N)$ converge to zero almost-everywhere on $R^N \setminus supp f$.

## ABSTRAK

Dalam makalah ini dikaji prinsip penempatan teritlak untuk kembangan spektrum bagi pengoperasi poliharmonik yang bertindih dengan kamiran Fourier berganda yang dijumlahkan ke atas domain yang sepadan dengan aras permukaan polinomial poliharmonik. Dibuktikan bahawa jumlah separa kamiran Fourier berganda bagi fungsi $f \in L_2(R^N)$ bertumpu kepada sifar hampir di mana-mana pada $R^N \setminus supp f$.

*Kata kunci*: Kamiran Fourier; kembangan spektrum; pengoperasi maksimum

## 1. Introduction

The well known classical localisation Theorem of Riemann states that if $2\pi-$periodic function $f$ is integrable over the interval $[-\pi, \pi]$, then the convergence or divergence of a one dimensional Fourier series at a given point depends only on the behaviour of the function $f \in L_1$ in an arbitrary small neighborhood of this point. The new concept of localisation principle has been introduced in Il'in (1968) and it is called generalised principle of localisation: the almost everywhere convergence of the spectral expansions of the given function to $0$ on the set, where expanding function vanishes. It is interesting to know whether the sufficient conditions for localisation can be weakened if we consider the general localisation problem. In this work we investigate the problems of the generalised localisation for the multiple Fourier integrals corresponding to the polyharmonic operator.

For any positive integer number *m*, we denote by $(-\Delta)^m$ a polyharmonic operator with the domain of definition $C_0^\infty(R^N)$, where $\Delta = \sum_{j=1}^{N} \frac{\partial^2}{\partial x_j^2}$ is the Laplace operator. As a symmetric operator the polyharmonic operator $(-\Delta)^m$ has at least one self adjoint extension in $L_2(R^N)$. Let *A* denote any self adjoint extension of the operator $(-\Delta)^m$ in $L_2(R^N)$. We consider the corresponding partial integrals of the multiple Fourier integrals:

$$E_\lambda f(x) = (2\pi)^{-\frac{N}{2}} \int_{|\xi|^{2m}<\lambda} \hat{f}(\xi) e^{i(x,\xi)} d\xi, \tag{1}$$

where $\hat{f}(\xi)$ is the Fourier transform of a function $f \in L_p(R^N)$, $1 \leq p \leq 2$.

Let $A$ denotes an unique self adjoint extension of the polyharmonic operator $\Delta^m$ in $L_2(R^N)$. Then the family $\{E_\lambda\}$ - be a decomposition of the identity of $A$ and the corresponding spectral expansion of any function $f \in L_2(R^N)$ coincides with (1) (for more details we refer the reader to Alimov *et al.* (1992)).

For the case of the Laplace operator $\Delta$ the generalised localisation principle (the almost everywhere convergence of the spectral expansions on $R^N \setminus supp(f)$) in classes $L_p(R^N)$ was investigated by many authors (Bastys 1983; Carleson 1966; Carbery & Soria 1988; 1997; Sjölin 1983). The problem of generalised localisation for spectral expansions of the Laplace operator is completely solved in classes $L_p(R^N)$, when $1 \leq p \leq 2$: if $p = 2$ then we have the generalised localisation and if $1 \leq p < 2$ then we do not. The interesting case of unit sphere was investigated in Anvarjon (2009). It should be noted that the problem of the generalised localisation for Riesz means of the Fourier-Laplace series on the critical line firstly investigated in this work.

In this paper we prove the following:

**Theorem 1.1.** *If $f \in L_2(R^N)$, and $f(x) = 0$ on an open set $\Omega \subset R^N$, then partial integrals*

$$E_\lambda f(x) = (2\pi)^{-\frac{N}{2}} \int_{|\xi|^{2m}<\lambda} \hat{f}(\xi) e^{i(x,\xi)} d\xi$$

*converge to $0$ almost everywhere on $\Omega$.*

In other words, the theorem asserts that for the multiple Fourier integrals the principle of generalised localisation holds in the class $L_2(R^N)$. We note that in case of Laplace operator a similar statement has been established in Carbery and Soria (1988). We originates from the latter paper, in which the authors obtained the relevant inequalities for maximal operator of the spectral decomposition of the Laplace operator. The maximal operator is defined by

$$E_* f(x) = \sup_{\lambda>0} |E_\lambda f(x)| \tag{2}$$

Theorem 1.1 is based on the following estimates of the maximal operator.

**Theorem 1.2.** *If $f \in L_2(R^N)$ and $f$ is supported in $\{|x| \geq 3\}$, then for any $r < 3$ there exists a constant $c_r$ such that*

$$\int_{|x|\leq r} \{E_* f(x)\}^2 dx \leq c_r \int_{|x|\geq 3} |f(x)|^2 dx. \tag{3}$$

The estimation of maximal operator can be applied in proofs concerning the almost everywhere convergence of spectral expansions. The almost everywhere convergence of spectral expansions were studied by many authors (Ashurov 1983; Anvarjon 2009; Carleson 1966; Sjölin 1983; Stein 1958); see for a review in Alimov (1970), and Zhizhiashvili and Topuriya (1979).

## 2. The Proof of Theorems

Let $f \in L_2(R^N)$, and $supp(f) \subset \{x \in R^N : |x| \geq 3\}$. We prove that for every $r < 3$ the following inequality holds:

$$\int_{|x|\leq r} \sup_{\lambda>0} |E_\lambda^{i\tau} f(x)|^2 dx \leq C_r \int_{|x|\geq 3} |f(x)|^2 dx, \tag{4}$$

with $\tau : -\infty < \tau < \infty$, here $E_\lambda^{i\tau} f(x) = (2\pi)^{-\frac{N}{2}} \int_{|\xi|^{2m}<\lambda} \left(1 - \frac{|\xi|^{2m}}{\lambda}\right)^{i\tau} \hat{f}(\xi) e^{i(x,\xi)} d\xi$. It should be noted here that the statement of the Theorem 1.2 is particular case of the latter inequality, when $\tau = 0$.

Let $\chi_b(t)$ be the characteristic function of the segment $[0,b]$ and $\phi(t)$ be a smooth function defined for $t \geq 0$, such that $\chi_{(3-r)/3}(t) \leq \phi(t) \leq \chi_{2(3-r)/3}(t)$. Then we define $\psi(x) = \phi(|x|) - \phi(2|x|)$ and $\psi_j(x) = \psi(\frac{x}{2^j})$ for $j = 1, 2, \ldots$. We obtain

$$\phi(|x|) + \sum_{j\geq 1} \psi_j(x) \equiv 1.$$

Let $\Theta_\lambda^{\tau,j} = \Theta_\lambda^{i\tau} \psi_j$. If $supp(f) \subset \{|x| \geq 3\}$, then for all $x : |x| \leq r$ we have

$$E_\lambda^{i\tau} f(x) = \Theta_\lambda^{i\tau} * f = \sum_{j\geq 1} \Theta_\lambda^{\tau,j} * f(x),$$

because $(\Theta_\lambda^{i\tau} \phi(|\cdot|) * f)(x) = 0$ if $|x| < r, r < 3$.

It is clear that to prove the inequality (4) it is sufficies to prove

$$\int_{R^N} \sup_{\lambda>0} |\Theta_\lambda^{\tau,j} g(x)|^2 dx \leq C\, 2^{-j} \int_{R^N} |g(x)|^2 dx, \quad \forall g(x) \in L_2(R^N).$$

By duality we can prove that the latter inequality can be established by investigating the Fourier transform of the "spectral function" $\Theta_t^{\tau,j}(x)$. Let $m_t^{\tau,j}(\xi) = (\hat{\Theta}_t^{\tau,j})(\xi)$. When $\tau = 0$ we use notation $m_t^j$, i.e. $m_t^{0,j}(\xi) = m_t^j(\xi)$. For $m_t^j(\xi)$ we have

**Lemma 2.1.** *For any $t > 0, \xi \in R^N, j \geq 1$ we have*

$$|m_t^j(\xi)| \leq \int_{||\xi|-t|<|y|2^{-j}} |\hat{\psi}(y)| dy.$$

This lemma is proved in Carbery and Soria (1988).

**Lemma 2.2.** *For any positive integer n there exists a constant C such that for any $t > 0, \xi \in R^N, j = 1,2,...$ we have the following estimate*

$$|m_t^{\alpha,j}(\xi)| \leq \frac{C}{(1+||\xi|-t|2^j)^n}.$$

**Proof**. We consider the following function:

$$m_t^{\alpha,j}(\xi) = \int_0^{\lambda} (1 - \frac{t}{\lambda})^{i\alpha} dm_t^j(\xi).$$

The integral we divide into two parts as follows

$$m_t^{\alpha,j}(\xi) = \int_0^{\lambda/2} (1 - \frac{t}{\lambda})^{i\alpha} dm_t^j(\xi) + \int_{\lambda/2}^{\lambda} (1 - \frac{t}{\lambda})^{i\alpha} dm_t^j(\xi) = I_1(\lambda) + I_2(\lambda).$$

Estimation of $I_2(\lambda)$. Using the following formula

$$\frac{d}{dt} m_t^j(\xi) = \frac{N}{t} m_t^j(\xi) + t^N \int_{|y|<2^j} \nabla \hat{\psi}(2^j \xi + yt) y dy,$$

we obtain :

$$|\frac{d}{dt} m_t^j(\xi)| \leq C 2^j \int_{||\xi|-t|2^j<|y|} (|\hat{\psi}(y)| + |\nabla \hat{\psi}(y)|(1+|y|)) dy.$$

Then for $I_2(\lambda)$ we have

$$|I_2(\lambda)| \leq \frac{c2^j}{\lambda} \int\limits_{\lambda/2}^{\lambda} \int\limits_{||\xi|-t|2^j<|y|} (|\hat{\psi}(y)| + |\nabla\hat{\psi}(y)|(1+|y|))dydt \leq$$

$$\leq \frac{c2^j}{\lambda} \int\limits_{\lambda/2}^{\lambda} \int\limits_{||\xi|-t|2^j}^{\infty} \frac{r^{N-1}}{(1+r)^l} drdt.$$

By changing the order of integration, and taking into account the relation
$\{(t,r): \lambda/2 \leq t \leq \lambda, ||\xi|-t|2^j \leq r \leq \infty\} \subset \{(t,r): ||\xi|-t|2^j \leq r \leq \infty, |\xi|-r \leq t \leq |\xi|+r\}$

$$|I_2(\lambda)| \leq \frac{c2^j}{\lambda} \int\limits_{||\xi|-t|2^j}^{\infty} \frac{r^{N-1}}{(1+r)^l} \int\limits_{|\xi|-r}^{|\xi|+r} dtdr = \frac{c2^j}{\lambda} \int\limits_{||\xi|-\lambda|2^j}^{\infty} \frac{r^N}{(1+r)^l} dr \leq$$

$$\leq \frac{c2^j}{\lambda} \frac{1}{(1+||\xi|-\lambda|2^j)^n}, n = l - N - 1.$$

Now we consider $I_1(\lambda)$. Using the integration by parts we obtain

$$I_1(\lambda) = \int\limits_0^{\lambda/2} (1-\frac{t}{\lambda})^{i\alpha} dm_t^j(\xi) = (1-\frac{t}{\lambda})^{i\alpha} m_t^j(\xi)\Big|_0^{\lambda/2} + \frac{i\alpha}{\lambda} \int\limits_0^{\lambda/2} (1-\frac{t}{\lambda})^{i\alpha}(1-\frac{t}{\lambda})^{-1} m_t^j(\xi) dt.$$

Using the inequality

$$|m_t^j(\xi)| \leq \int\limits_{||\xi|-t|\leq|y|} |\hat{\psi}(y)| dy,$$

we get

$$|I_1(\lambda)| \leq m_{\lambda/2}^j(\xi) + \frac{2|\alpha|}{\lambda} \int\limits_0^{\lambda/2} \int\limits_{||\xi|-t|\leq|y|} |\hat{\psi}(y)| dydt.$$

Similarly as in previous case we have

$$\int\limits_0^{\lambda/2} \int\limits_{||\xi|-t|\leq|y|} |\hat{\psi}(y)| dydt \leq C \int\limits_0^{\lambda/2} \int\limits_{||\xi|-t|}^{\infty} \frac{r^{N-1}}{(1+r)^l} drdt \leq$$

$$\leq C \frac{1}{(1+\varepsilon_0 ||\xi|-\lambda|2^j)^n}.$$

Consequently for $I_1(\lambda)$ we have

$$|I_1(\lambda)| \leq \frac{C}{(1+||\xi|-\lambda|2^j)^n}.$$

Finally we obtain

$$|m_\lambda^{\alpha,j}(\xi)| \leq \frac{C}{(1+||\xi|-\lambda|2^j)^n}.$$

The lemma is proved. □

For the derivative of the function $m_t^{\alpha,j}(\xi)$ we obtain

**Lemma 2.3.** *For any positive integer n there exists a constant C such that for any $t > 1$ and $\xi \in R^N$, $j = 1, 2,...$ we have the following estimates*

$$|\frac{d}{dt}m_t^{\alpha,j}(\xi)| \leq \frac{C2^j}{(1+\varepsilon_0||\xi|-t|2^j)^n}.$$

Note from (4) by putting $t = 0$ we have the statement of Theorem 1.2.

Proof of Theorem 1.1. Let $f \in L_p(R^N)$ and $f(x) = 0$ on an open set $\Omega \subset R^N$. We have to prove that $E_\lambda f(x) \to 0$ almost everywhere on $\Omega$, or if $x_0 \in \Omega$ is an arbitrary point, then on the ball $B_{r_0}(x_0) \subset \Omega$. Therefore without loss of generality we assume that $supp(f) \subset R^N \setminus B_{r_0}(x_0)$. Due to the invariance of $E_\lambda f(x)$ with respect to the transaction and dilation, we can reduce the problem to the consideration of the functions with $supp(f) \subset \{|x| \geq 3\}$.

Thus we have to prove that for any function $f \in L_2(R^N)$ with $supp(f) \subset \{|x| \geq 3\}$, one has $E_\lambda f(x) \to 0$, almost everywhere on $\{|x| < r\}, r < 3$. But this is a consequence of inequality for $E_\lambda f$ in Theorem 1.2, because the latter inequality allows us to state that the set $\{x : E_\lambda f(x)$ does not converge to $f(x)\}$ has a measure 0 (for more details see Stein and Weiss (1971)). Theorem 1.1 is completely proven. □

**Acknowledgement**

This research has been supported by Universiti Putra Malaysia under Research University Grant (RUGS): 05-01-09-0674RU.
    The first author wishes to thank Prof. Ashurov R.R. for his kind advice and stimulating discussions.

**References**

Alimov Sh.A., Ashurov R.R. & Pulatov A.K. 1992. Multiple Fourier series and Fourier integrals. In: Khavin V.P. & Nikol'skii N.K.(Eds.). *Commutative Harmonic Analysis IV*. New York: Springer-Verlag.

Alimov Sh.A. 1970. Summability almost everywhere of Fourier series in Lp with respect to eigenfunctions. *Journal of Differential Equations* **6**(1): 164-171.

Anvarjon Ahmedov. 2009. The principle of general localization on unit sphere. *Journal of Mathematical Analysis and Applications* **356**(1): 310-321.

Ashurov R.R.1983. Summability almost everywhere of Fourier series in Lp with respect to eigenfunctions. *Mathematical Notes of the Academy of Sciences of the USSR* **34**(6): 913-916

Bastys A.J. 1983. Generalized localization of Fourier series with respect to the eigenfunctions of the Laplace operator in the classes Lp. *Litovskii Matematicheskii Sbornik* **31**(3): 387-405.

Carleson L. 1966. On convergence and growth of partial sums of Fourier series. *Acta Mathematica* **116**: 135-157.

Carbery A. & Soria F. 1988. Almost everywhere convergence of Fourier integrals for functions in Sobolev spaces, and an $L_2$-localization principle. *Revista Mat. Iberoamericana* **4**(2): 319–337.

Carbery A. & Soria F. 1997. Pointwise Fourier inversion and localisation in $R^n$. *Journal of Fourier Analysis and Applications* **3** (Special issue): 847-858.

Il'in V.A. 1968. Localization and convergence problems for Fourier series by fundamental function systems of the Laplace operator. *Russian Math. Surveys* **23**(2): 59-116.

Sjölin P. 1983. Regularity and integrability of spherical means. *Monatsh. Math.* **96**(4): 277-291.

Stein E.M. & Weiss G. 1971. *Introduction to Fourier Analysis on Euclidean Spaces*. Princeton, NJ: Princeton Univ. Press.

Stein E.M. 1958. Localization and summability of multiple Fourier series. *Acta Mathematica* **100**: 93-147.

Zhizhiashvili L.V. & Topuriya S.B. 1979. Fourier-Laplace series on a sphere. *Journal of Soviet Mathematics* **12**(6): 682-714.

*Department of Process and Food Engineering*
*Faculty of Engineering*
*Universiti Putra Malaysia*
*43400 UPM Serdang*
*Selangor DE, MALAYSIA*
*E-mail: anvar@eng.upm.edu.my*